# Stratégies de contrôle pour les éoliennes flottantes : état de l'art et perspectives

Flavie DIDIER, Salah LAGHROUCHE, Daniel DEPERNET
FEMTO-ST Institute, Univ Bourgogne Franche-Comté, UTBM, CNRS, Belfort, France
*flavie.didier@utbm.fr*

**RESUME** – Le secteur de l'éolien flottant possède un potentiel énergétique conséquent. Cependant, la minimisation des mouvements de la structure sous l'effet conjugué du vent et des vagues tout en assurant une extraction maximale de la puissance sur une large plage de fonctionnement constitue l'un des principaux défis pour le pilotage de ces éoliennes. Cet article propose une revue des méthodes de contrôle issues de la littérature pour les éoliennes flottantes. Les limitations de ces contrôleurs sont discutées, avant d'introduire une présentation de plusieurs méthodes prometteuses basées sur des données. En particulier, ce papier met l'accent sur les techniques d'intelligence artificielle associées aux méthodes de contrôle à base de données. Enfin, le projet CREATIF traitant de la simulation en temps réel d'éoliennes flottantes et de leurs contrôles intelligents est présenté.

**ABSTRACT** – The floating wind turbines sector has great energy potential. However, minimizing the movement of the structure under the combined effect of wind and waves while ensuring maximum power extraction over a wide operating range is one of the main challenges for the control of these turbines. This paper presents a review of control methods for floating wind turbines from the recent literature. The limitations of these controllers are discussed, before introducing a presentation of several promising data-based methods. In particular, this paper focuses on artificial intelligence techniques associated with data-based control methods. Finally, the CREATIF project dealing with real-time simulation of floating wind turbines and their intelligent controls is presented.

**MOTS-CLES** – éolienne flottante, énergie marine renouvelable, contrôle robuste et adaptatif, intelligence artificielle

## 1. Introduction

Dans un souci de réduction des émissions de gaz à effet de serre et ainsi de réduction de la consommation en ressources fossiles, la transition énergétique en faveur des énergies renouvelables se révèle être une alternative prometteuse. En particulier la technologie éolienne flottante qui connaît un véritable essor depuis quelques années.

En effet, l'éolien flottant est une filière émergente à fort potentiel permettant de placer des parcs éoliens plus loin des côtes, les rendant ainsi moins visibles et déployables en plus grand nombre. Depuis l'installation en Europe de la première éolienne flottante en 2009, plusieurs projets de parcs ont vu le jour, avec notamment le premier parc commercial « Hywind » (2017) en Écosse. En tout, la puissance installée en Europe d'ici 2024 est estimée à 312,8 MW. La France se positionne sur ce secteur avec une machine pilote (« Floatgen ») de 2 MW de capacité en exploitation depuis 2018, ainsi que quatre parcs (115 MW de capacité) prévus pour 2023. La Programmation pluriannuelle de l'énergie prévoit des appels d'offres pour la construction de parcs en Méditerranée et en Bretagne. Le potentiel de production électrique en France métropolitaine pour l'éolien flottant est estimé à 155 GW [1].

Aujourd'hui, cette filière doit faire face à de multiples défis. L'un des principaux défis est le développement de systèmes de contrôle robustes, adaptatifs, et intelligents permettant la minimisation des mouvements de la structure tout en assurant une extraction maximale de la puissance sur l'ensemble des zones de fonctionnement de l'éolienne.

Ce papier propose donc une revue non exhaustive des principales méthodes de la littérature récente pour le pilotage des éoliennes flottantes en mettant l'accent sur les défis et verrous inhérents à ce type de système. De plus, il propose un exposé des perspectives et challenges pour surmonter ceux-ci. L'article est organisé comme suit. Dans la section 2, une description du système éolien flottant, ainsi que les objectifs liés au système de contrôle sont introduits. Une classification des méthodes de pilotage de ces éoliennes est proposée dans la section 3. Les méthodes ont été organisées selon le type de contrôle, conventionnel ou avancé. Dans la section 4, l'approche basée sur les données, avec notamment l'utilisation de méthodes d'intelligence artificielle, est étudiée comme une perspective intéressante de surpasser les limitations des méthodes actuelles. Une présentation du projet CREATIF et de ses ambitions sont exposées dans la section 5.

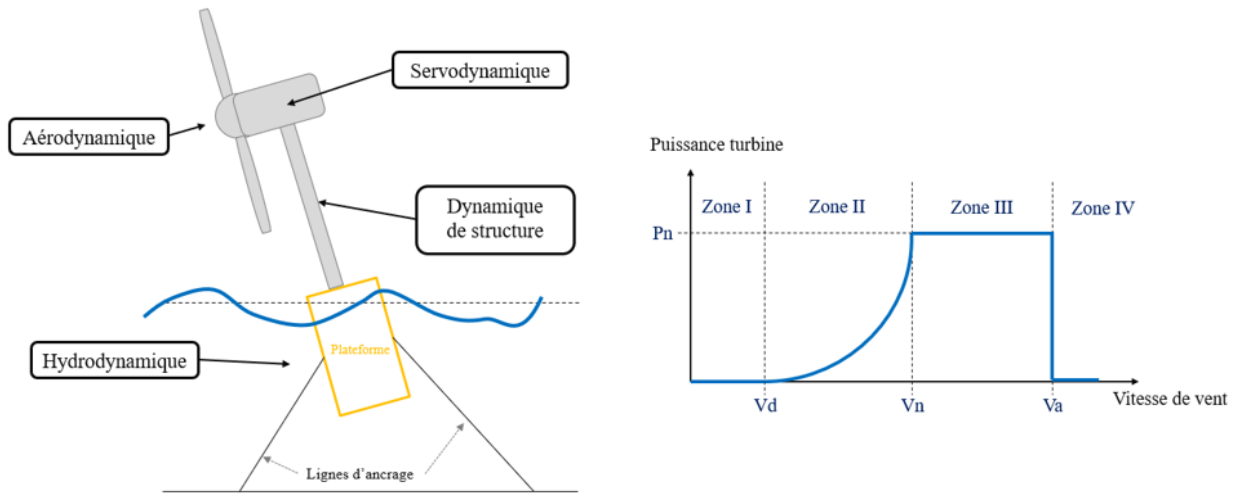

**Figure 1 : a. Éolienne flottante - b. Zones de fonctionnement**

## 2. Description du système

### 2.1 L'éolienne flottante

Contrairement à l'éolien en mer, la technologie de l'éolien flottant permet de s'affranchir de la contrainte de profondeur ne pouvant dépasser 60 mètres. De plus, les éoliennes peuvent être situées dans des zones avec un meilleur gisement de vent favorisant ainsi une meilleure production d'électricité. La structure d'une éolienne flottante est composée d'une part de l'éolienne, pour capter l'énergie du vent, et de l'autre de la plateforme flottante accompagnée des lignes d'ancrage, pour assurer le support de l'éolienne. La plateforme flottante sur laquelle repose l'éolienne est inspirée des technologies utilisées dans les secteurs des plateformes pétrolières et gazières. Plusieurs types de flotteurs cohabitent sur le marché pour l'éolienne flottante, à savoir, barge, immergé à câbles tendus (TLP), à espar, et semi-submersible [2].

Les principales dynamiques (figure 1.a) s'appliquant sur la structure sont les charges aérodynamiques liées au vent, les charges hydrodynamiques liées aux vagues, les charges de structure s'appliquant sur l'éolienne et les charges liées aux systèmes électriques au sein de la nacelle (servo-dynamique). Contrairement à l'éolien fixe, la plateforme flottante introduit donc des charges hydrodynamiques supplémentaires en raison de la houle. De plus, un couplage entre les charges aérodynamiques et hydrodynamiques existe, faisant de l'éolien flottant un système complexe.

### 2.2 Structure et objectif de commande

La plage d'opération de l'éolienne flottante est divisée en trois zones de fonctionnement selon la vitesse de vent entrant (figure 1.b). Dans la zone I, la vitesse de vent est inférieure à la vitesse de démarrage de l'éolienne, aucune puissance électrique n'est produite. Dans la zone II, la vitesse de vent est comprise entre la vitesse de démarrage et la vitesse nominale, l'objectif de contrôle est alors d'extraire le maximum de puissance à partir du vent, tout en stabilisant l'éolienne afin de préserver sa structure. Ces objectifs sont réalisés par la commande du couple de la génératrice et de l'angle des pales. Dans la zone III, la vitesse de vent est supérieure à la vitesse nominale, l'objectif est alors de réguler la puissance générée, en la limitant à la puissance nominale, par variation de l'angle des pales. Lorsque la vitesse de vent dépasse une vitesse limite, les pales de l'éolienne sont mises en bateau pour éviter d'endommager la structure. Un mécanisme de contrôle adéquat est nécessaire pour atteindre les objectifs en termes de puissance et faire face aux mouvements de la plateforme.

## 3. Etat de l'art des méthodes de commande pour l'éolien flottant

Aujourd'hui, peu de solutions de contrôle viables existent pour l'éolien flottant, en effet, le problème de sa commande est encore récent dans la communauté automaticienne. Cette section introduit les méthodes de contrôle conventionnel empruntées aux systèmes d'éoliennes fixes, puis s'étend aux méthodes de contrôle avancé pour l'éolien flottant.

### 3.1 Méthodes de contrôle conventionnel pour l'éolien flottant

Les contrôleurs conventionnels d'éolienne flottante sont des mécanismes de contrôle simples et faciles à concevoir basés sur le principe « single-input single-output » (système SISO). Ils ne requièrent pas de modèle analytique du système, seulement des boucles fermées avec mesure du signal de sortie. Ainsi, des boucles de contrôle indépendantes sont appliquées en parallèle pour atteindre les objectifs de contrôle selon la zone de fonctionnement.

Dans le cas de l'éolien flottant, afin de minimiser le mouvement de tangage de la plateforme, la fréquence de réglage du pas des pales est gardée plus basse que la fréquence de résonnance de la plateforme [3]. Pour ce faire, une commande par séquencement de gains (GSPI) pour la zone III a été testée. Un autre contrôleur GSPI permettant limiter le phénomène d'amortissement négatif a été utilisé en [4] ajustant cette fois-ci la vitesse du rotor en fonction de l'activité de l'angle de pale. Ces contrôleurs avec des gains « désaccordés » par rapport à l'éolien fixe se sont avérés être les contrôleurs les plus appropriés, car ils réduisent l'activité des pales et résolvent le problème du tangage de la plateforme. Cette configuration est aujourd'hui considérée comme la commande de référence des éoliennes flottantes, et elle est utilisée pour tester de nouveaux contrôleurs. D'ailleurs, ce contrôleur de référence a été analysé pour différents types de plateformes [5].

L'instabilité de la plateforme a également été étudiée en utilisant la vitesse de tangage de la plateforme comme entrée pour réguler la vitesse nominale du générateur en zone III. Dans [6] une stratégie de contrôle basée sur l'estimation de la vitesse de vent est proposée pour supprimer l'amortissement négatif. L'efficacité de cette méthode est gouvernée principalement par la qualité de cette estimation.

Les méthodes conventionnelles empruntées à l'éolien fixe ne peuvent constituer qu'un point de départ. En effet, elles sont limitées en termes de performance et de gestion du couplage des boucles de contrôle ainsi que des perturbations (vent et vague), et ne sont donc pas les plus adaptées pour un système éolien flottant fortement couplé et multi-objectif. Ainsi, plusieurs travaux proposent d'utiliser des approches de synthèse de commandes basées sur l'utilisation de modèles « multiple-input multiple-output » (systèmes MIMO) qui prennent en compte les couplages des sorties à piloter.

## 3.2 Méthodes de contrôle avancé pour l'éolien flottant

En général, la réalisation de ces méthodes de contrôle repose sur des modèles analytiques des dynamiques du système. À ce jour, les travaux les plus aboutis et exploitables pour la commande sont ceux de Betti [7], et de Sandner [8]. Ces approches de modélisations non linéaires sont toutes deux basées sur des hypothèses restrictives afin d'obtenir des modèles non linéaires simplifiés. En général, pour synthétiser leurs commandes les chercheurs passent par une linéarisation autour d'un point de fonctionnement du modèle.

Les principaux contrôleurs proposés dans la littérature récente sont, la commande linéaire quadratique (LQR), la commande H-infini, la commande de systèmes linéaires à paramètres variants (LPV) et la commande prédictive (MPC). Tous ont été développés en se basant sur une linéarisation autour d'un ou plusieurs points de fonctionnement. Les contrôleurs LQR sont notamment utilisés pour diminuer le mouvement de la plateforme. Des contrôleurs collectifs du pas des pales (LQR-CBPC) et individuel du pas des pales (LQR-IBPC) ont été conçus [9], LQR-IBPC permet de meilleures performances en termes de vitesse du rotor et de régulation de puissance mais le mouvement de tangage de la plateforme reste important. Une amélioration possible est le contrôleur LQR-IBPC avec commande par adaptation à la perturbation, grâce à sa capacité d'estimation du vent il minimise son action sur le système, assurant ainsi une meilleure régulation en vitesse et puissance [10]. De même, un contrôleur H-infini-CPB a été développé pour gérer le mouvement de la plateforme, et aussi pour la régulation de puissance et la réduction des charges associées [11]. Un contrôleur avec séquencement de gain (GS) et stratégie de sortie H-infini a notamment permis une amélioration en termes de charge sur le mât et de régulation de vitesse du rotor [12]. Dans [13], un LPV par retour d'état et un contrôleur LQR par retour de sortie basés tous deux sur GS sont proposés pour réguler la puissance et minimiser la charge structurelle. Ils ont montré de meilleures performances comparées au régulateur de référence. Un contrôleur CBP-LPV a aussi été proposé dans [14].

D'autres méthodes de synthèse basées sur les systèmes non linéaires comme la stratégie de linéarisation entrées-sorties (IOFL) et la commande par mode glissant (SMC) ont aussi été utilisées pour réguler la vitesse de la génératrice tout en analysant les effets de perturbations entrantes sur le mouvement de la plateforme [15]. Une amélioration par rapport au contrôleur de référence en termes de régulation de la vitesse de la génératrice a été constatée pour SMC, le mouvement de la plateforme reste quant à lui similaire. En revanche, IOFL entraîne une augmentation du mouvement de la plateforme. Néanmoins, ces méthodes ont été synthétisées en utilisant un modèle non linéaire simplifié en supposant une connaissance parfaite de ses paramètres.

Un contrôleur MPC-CPB basé sur un modèle linéaire opérant pour réguler la puissance à une valeur constante et minimiser les charges structurelles a démontré de meilleurs résultats que le contrôleur PI en termes de régularisation de la puissance et de réduction du tangage de la plateforme [16]. Parmi les approches MPC basées sur des modèles non linéaires, on peut citer la commande prédictive non linéaire (NMPC-CBP) basé sur le modèle simplifié de Sandner [17]. Ce contrôleur utilise l'estimation du vent et des vagues, et obtient des résultats satisfaisants pour la régulation de la vitesse et de la puissance ainsi que la réduction de charge sur les pales. Cependant, il demande plus de ressources computationnelles. Un contrôleur NMPC (IBP) a par la suite été proposé permettant de réduire d'avantage les charges sur les pales et sur le rotor [18]. L'applicabilité de ces contrôleurs pose question quant à leur capacité à piloter le système en temps réel.

### 3.3 Limites

Une représentation fidèle de l'éolienne flottante passe par une modélisation multi-physiques complexe non linéaire. De plus, toutes les dynamiques ne peuvent être modélisées parfaitement (incertitudes et dynamiques non modélisées) impactant alors la robustesse des contrôleurs. Paradoxalement, plus le modèle est précis, plus la conception du système de commande est coûteuse et ainsi une réduction du modèle ou du contrôleur doit alors être réalisée pour être utilisable en pratique. La plupart des contrôleurs pour l'éolien flottant impliquent donc la linéarisation autour d'un point de fonctionnement du modèle, où s'éloigner de ce point peut entraîner une dégradation des performances.

## 4. Approche par les données et intelligence artificielle

Avec l'avancée des technologies en matière de capteurs, d'outils de simulation, ou encore de stockage, beaucoup de données peuvent être générées et stockées. Ces données, en ligne ou hors ligne, peuvent alors être utilisées pour concevoir directement le contrôleur. La question du contrôle pour l'éolien flottant étant relativement récente, cette section a pour intention de souligner les éventuelles approches à explorer.

### 4.1 Méthodes de contrôle basées sur les données

Les stratégies de contrôle pilotées par les données sont toutes les théories et méthodes de contrôle dans lesquelles le contrôleur est conçu en utilisant directement les données d'entrée/sortie du système contrôlé ou les connaissances issues du traitement des données, mais sans aucune information explicite du modèle mathématique du processus contrôlé [19]. La stabilité, convergence et robustesse peuvent être garanties par analyse mathématique rigoureuse.

Il existe une multitude de méthodes de contrôle pilotées par les données [19]. Une classification de celles-ci peut être faite selon le type de données utilisées, en ligne, hors ligne ou hybride, pour le contrôle de systèmes commandés linéaires ou non linéaires. Les méthodes basées sur les données mesurées en ligne sont notamment, l'approximation stochastique des perturbations simultanées (SPSA) où la structure du contrôleur est inconnue et elle repose sur une approche par identification des paramètres du contrôleur. De même, la méthode de contrôle non falsifié repose sur la même approche mais nécessite une structure de contrôleur fixée et connue pour le contrôleur. Un contrôleur non falsifié basé sur les données est développé pour le positionnement dynamique des systèmes marins [20]. La commande adaptative sans modèle (MFAC) [21] est basée sur le principe de linéarisation dynamique et ne nécessite pas de structure fixée pour le contrôleur. Celui-ci étant construit à partir des données entrées et sorties mesurées. Ces trois méthodes de commandes sont adaptatives et particulièrement adaptées aux systèmes linéaires.

Parmi les méthodes mobilisant des données hors ligne, autre que la commande PID, on peut aussi citer la méthode par réglage itératif de la rétroaction (IFT), la méthode par renforcement de la référence virtuelle par rétroaction (VRFT) et la méthode par réglage basée sur la corrélation (CbT). Elles requièrent une structure fixée du contrôleur et utilisent une approche par identification des paramètres.

Enfin, il existe des méthodes basées à la fois sur des données mesurées en ligne et hors ligne (hybrides). Parmi elles, la méthode de contrôle par apprentissage itératif (ILC) sans structure fixée du contrôleur au préalable, et la méthode adaptative d'apprentissage paresseux (LL) basée sur la linéarisation dynamique. Toutes les deux sont adaptées aux systèmes non linéaires.

Les méthodes les plus intéressantes pour l'éolien flottant sont les méthodes de contrôle adaptées aux systèmes non linéaires, comme MFAC, SPSA, contrôleur non falsifié ou LL. Les principales limitations sont le manque de méthodologie pour le choix de la structure fixée du contrôleur (SPSA et non falsifié), le manque d'analyse de la stabilité (SPSA, MFAC, LL) et un coût computationnel pouvant être élevé (SPSA, LL).

### 4.2 Méthodes de contrôle basées sur l'intelligence artificielle :

L'augmentation des capacités de calcul et de traitement de gros volume de données est une opportunité pour la mise en place de contrôleurs reposant sur des techniques dites d'intelligence artificielle (IA). En matière de contrôle, les techniques d'IA sont essentiellement basées sur des techniques de régression ou d'optimisation [22]. Pour la régression, on retrouve principalement les systèmes experts, la logique floue et l'apprentissage automatique (« Machine Learning »). Alors que pour l'optimisation, l'apprentissage machine et les méthodes métaheuristiques sont à privilégier. Aujourd'hui, les systèmes experts et la logique floue sont mis à l'écart aux profits des méthodes plus puissantes de métaheuristiques et d'apprentissage automatique. En particulier, l'apprentissage par renforcement profond connaît un essor d'attractivité pour son application dans le domaine du contrôle [23].

En effet, l'apprentissage par renforcement est particulièrement adapté aux systèmes avec peu de connaissances ou avec un modèle difficile à formuler. Le principe repose sur l'interaction entre un agent en mesure de prendre des décisions dans l'espace d'action, et un environnement caractérisé par l'espace d'état. La qualité de l'action prise par l'agent est évaluée à l'aide d'une stratégie de récompense. Une politique définit le comportement de l'agent selon l'état et une fonction valeur attribue à chaque état sa valeur (basée sur les récompenses attendues) [24]. Deux approches coexistent, celle basée sur la valeur et celle basée sur la politique. On parle d'apprentissage par renforcement profond, lorsque l'approximation de la fonction de valeur ou de la politique se fait par l'intermédiaire d'une technique d'apprentissage supervisé, telles que les réseaux de neurones. L'apprentissage par renforcement profond permet de s'affranchir de la discrétisation de l'espace d'état et possède une bonne propriété de généralisation face aux nouvelles situations.

Les principaux algorithmes basés sur la valeur sont dérivés de l'algorithme d'apprentissage-Q (Q-learning), on peut citer l'apprentissage-Q ajusté ou encore l'apprentissage-Q profond (DQN). Un contrôleur MPPT pour un système de conversion d'énergie éolienne à vitesse variable basé sur un algorithme de Q-learning est proposé dans [25]. L'apprentissage d'une cartographie optimale des états aux actions est réalisé en ligne grâce à la mise à jour des valeurs d'action selon les récompenses futures. Cette approche n'exige aucune connaissance des paramètres de l'éolienne, puisque l'agent est capable d'apprendre par interaction directe avec l'environnement. Dans [26] un contrôleur MPPT pour l'éolien a également été étudié cette fois-ci avec un algorithme DQN, afin d'apprendre la relation optimale entre vitesse du rotor et puissance électrique sortante.

Toutefois, ces algorithmes ne sont pas les plus adaptés pour gérer les problèmes avec un espace d'action continu et peuvent devenir instables. En revanche, les algorithmes basés sur le gradient et en particulier par acteur-critique permettent l'extension aux espaces d'action continus. Les principaux algorithmes sont le gradient de politique déterministe profond (DDPG), l'avantage acteur-critique (A2C) ou encore acteur-critique doux (SAC). De plus, ceux-ci possèdent une plus forte garantie de convergence et un apprentissage plus rapide que les algorithmes basés sur la valeur.

L'apprentissage par renforcement profond peut être l'outil le plus approprié pour les problèmes de contrôle, car il permet la synthèse de contrôleurs adaptatifs et une faculté de généralisation aux nouvelles situations rencontrées. Contrairement aux algorithmes de DQN qui requièrent la discrétisation des actions, les algorithmes de DDPG peuvent gérer les espaces de hautes dimensions. Il semble donc intéressant d'étudier plus en avant ces algorithmes qui peuvent assurer un pilotage des éoliennes flottantes avec un bon niveau de performance.

## 5. Le projet CREATIF

Dans le cadre du projet CREATIF, on s'intéresse au contrôle opérationnel d'une éolienne flottante sur l'ensemble des zones de fonctionnement II et III, avec une attention particulière portée sur la transition entre ces deux zones. L'éolienne flottante étudiée est une éolienne à 3 pales de modèle DTU 10 MW. Elle repose sur une plateforme de type espar et elle est installée dans la région du golf du Maine.

L'ambition de ce projet est non seulement la synthèse de lois de commande pour l'éolienne flottante mais également leurs mises en œuvre. Pour ce faire, dans un premier temps les contrôleurs seront synthétisés et testés à l'aide du logiciel de simulation en temps réel OpenFAST développé par NREL. Par la suite des essais seront entrepris sur un banc de type « Power hardware in the loop » (PHIL) mis en place à FEMTO-ST, Belfort, permettant la simulation intégrée et en temps réel du système complet avec la dynamique de l'éolienne flottante simulée numériquement par OpenFAST, et la partie chaîne de conversion de l'énergie utilisant des dispositifs de puissance réels d'échelle réduite. Enfin, un dispositif évoluant dans un bassin à vagues situé au LHEEA (ECN Nantes) émulant le fonctionnement d'une éolienne flottante permettra la réalisation de test pour la mise en œuvre du contrôleur développé. Les principaux objectifs de contrôle sont de maximiser la production de puissance (avec prise en compte de la qualité de l'énergie) et de préserver la structure en réduisant les dommages liés à la fatigue.

Afin de synthétiser le contrôleur, le choix de l'approche est crucial. En effet, on souhaite un contrôleur adaptatif et robuste capable d'opérer sur l'ensemble des zones de fonctionnement et implantable en temps réel. Or, la plupart des méthodes actuelles de contrôle avancé sont limitées à la linéarisation autour d'une multitude de points de fonctionnement, rendant la conception du contrôleur fastidieuse et valide sur une plage de fonctionnement restreinte. C'est pourquoi, l'utilisation de l'apprentissage par renforcement profond constitue une alternative prometteuse pour le contrôle de systèmes complexes et incertains tels que les éoliennes flottantes.

## Références